\documentstyle[twoside,amssymb,12pt]{article}
\def\C{{\Bbb C}}

\def\N{{\Bbb N}}
\def\Q{{\Bbb Q}}
\def\R{{\Bbb R}}
\def\P{{\Bbb P}}
\def\Z{{\Bbb Z}}

\newtheorem{prop}{Proposition}[section]
\newtheorem{dfn}[prop]{Definition}
\newtheorem{theo}[prop]{Theorem}
\newtheorem{conj}[prop]{Conjecture}
\newtheorem{rem}[prop]{Remark}
\newtheorem{coro}[prop]{Corollary}
\newtheorem{lem}[prop]{Lemma}
\newtheorem{exam}[prop]{Example}

\title{\sc Non-Archimedian integrals and stringy Euler numbers of 
log terminal pairs} 

\author{{\sc Victor V. Batyrev} \\
\small  {\em Mathematisches Institut, Universit\"at T\"ubingen}   \\
\small  {\em Auf der Morgenstelle 10,  72076  T\"ubingen, Germany}  \\
\small  {\em  e-mail: batyrev@bastau.mathematik.uni-tuebingen.de} \\
}
 
\begin{document}

\date{}

\maketitle

\begin{abstract}
Using non-Archimedian integration  over spaces of arcs of
algebraic varieties,  we define stringy  Euler numbers associated with
arbitrary Kawamata log terminal pairs.   There is a natural  Kawamata
log terminal pair corresponding  to  an  algebraic variety  $V$
having a  regular action of a finite group $G$. In this  situation
we  show  that  the stringy Euler number of this pair coincides with
the physicists' orbifold Euler number defined by the Dixon-Harvey-Vafa-Witten
formula.   As an application, we   prove a conjecture of Miles Reid
on the Euler numbers of  crepant desingularizations
of Gorenstein quotient singularities.
\end{abstract}

\thispagestyle{empty}

%\newpage

\section{Introduction}

Let $X$ be a normal irreducible algebraic variety of dimension $n$ over 
$\C$, 
$Z_{n-1}(X)$ the 
group of Weil divisors on $X$,  ${\rm Div }(X) \subset Z_{n-1}(X) $ 
the subgroup of Cartier divisors on 
$X$, $Z_{n-1}(X) \otimes \Q$ 
the group of Weil divisors on $X$ with coefficients in $\Q$, 
$K_X \in Z_{n-1}(X)$ a  
canonical divisor of  $X$. 

Recall several definitions from the Minimal Model Program 
\cite{K1,K2,KMM} (see also 
\cite{kollar1,kollar2}):

\begin{dfn}
{\rm Let  $\Delta_X  \in  Z_{n-1}(X) \otimes 
\Q$ be a $\Q$-divisor on a normal irreducible algebraic 
variety $X$. A resolution of singularities  $\rho \,: \, Y  \rightarrow  X$
is called a {\bf log resolution of $(X, \Delta_X)$} if 
the union of the $\rho$-birational transform $\rho^{-1}(\Delta_X)$ of 
$\Delta_X$ with the exceptional 
locus of $\rho$ is a divisor $D$ consisting of  smooth irreducible 
components $D_1, \ldots, D_m$ having  only normal crossings. }
\end{dfn}

\begin{dfn} 
{\rm Let  $\rho \,: \, Y  \rightarrow  X$ be a log resolution of a 
pair $(X, \Delta_X)$. We assume that $K_X + \Delta_X$ is a $\Q$-Cartier 
divisor and write 
\[ K_{Y} = \rho^* (K_X + \Delta_X) + \sum_{i=1}^m a(D_i, \Delta_X) D_i,  \] 
where $D_i$ runs through 
all irreducible components of $D$ and $a(D_i, \Delta_X) = -d_j$ 
if $D_i$ is a $\rho$-birational transform of an irreducible component 
$\Delta_j$ of $Supp\, \Delta_X$.  Then the  number rational 
number $a(D_i, \Delta_X)$ (resp. 
$a_l(D_i, \Delta_X):=  a(D_i, \Delta_X) +1$) is called
the {\bf discrepancy} (resp. {\bf log discrepancy}) of $D_i$.} 
\end{dfn}

\begin{dfn} 
{\rm A pair  $(X, \Delta_X)$ 
is called 
{\bf  Kawamata log terminal}  
 if the following conditions are satisfied:   

 (i) $ \Delta_X  = d_1 \Delta_1 + \cdots + d_k \Delta_k$, where 
$\Delta_1, \ldots,\Delta_k$ are  irreducible Weil divisors  
and $d_i < 1$ for all $i \in \{ 1, \ldots, k \}$;  

(ii) $K_X + \Delta_X$ is a $\Q$-Cartier divisor; 

(iii)  for any log   resolution of singularities   
$\rho \,: \, Y  \rightarrow  X$, we have  
$a_l(D_i,\Delta_X) > 0 $ for all $i \in \{1, \ldots, m \}$,  }
\end{dfn}

Now we introduce a new invariant of Kawamata log terminal pairs:

\begin{dfn}
{\rm Let $(X, \Delta_X)$ be a Kawamata 
log terminal pair, $\rho\, : \, Y \rightarrow 
X$ a log resolution of singularities  as above. We 
put   $I = \{1, \ldots, m\}$ and set  for  any subset $J \subset I$  
\[ D_J := \left\{ \begin{array}{ll}   \bigcap_{ j \in J} D_j & \mbox{\rm 
if $J \neq \emptyset$} 
\\
Y  & \mbox{\rm 
if $J =  \emptyset$} \end{array} \right. , \;\;\; 
D_J^{\circ} := D_J \setminus \bigcup_{ j \in (I \setminus J)} D_j, $$
$$e(D_J^{\circ}):= (\mbox{\rm topological Euler number of  $D_J^{\circ}$}). $$
 We call the rational number  
\[ e_{\rm st}(X,\Delta_X) := \sum_{J \subset I} e(D_J^{\circ}) \prod_{j \in J} 
{a_l(D_j,\Delta_X)}^{-1} \]
the {\bf stringy  Euler number} of the Kawamata 
log terminal pair $(X, \Delta_X)$
(in the above formula, we assume   $\prod_{j \in J} = 1$ if $J = \emptyset$).} 
\label{ph}
\end{dfn}

Using non-Archimedian integrals, 
we show that the stringy Euler number 
 $e_{\rm st}(X, \Delta_X)$ is well-defined: 

\begin{theo} 
In the above definition,  
$e_{\rm st}(X, \Delta_X)$ does not depend on the choice of a 
log resolution  $\rho\, : \, Y \rightarrow X$. 
\label{indep}
\end{theo} 

We expect that the stringy Euler numbers have  the 
following natural connections with log flips in dimension $3$ 
(see \cite{Sh1,Sh2}): 

\begin{conj} 
Let $X$ be a normal $3$-dimensional variety and $\Delta_X$ is 
an effective $\Q$-divisor such that $(X, \Delta)$ is Kawamata log terminal,  
and  $\varphi\,:\,  (X, \Delta_X)  \dasharrow (X^+, \Delta_{X^+})$ a 
log $(K_X + \Delta_X)$-flip.  Then one has the following inequality: 
\[ e_{\rm st}(X, \Delta_X) > e_{\rm st}(X^+, \Delta_{X^+}). \]
 \end{conj} 

\begin{rem} 
{\rm In \ref{lflip2} we show that the above conjecture is true for 
toric log flips in arbitrary dimension 
$n$.}    
\end{rem} 

Recall now a  definition from  the string theory \cite{DHVW} (see also 
\cite{R}): 

\begin{dfn}
{\rm Let $V$ be a  smooth complex algebraic  variety together  with a regular 
action of a finite group $G$: $G \times V \rightarrow V$. For any element 
$g \in G$ we set 
$$V^g := \{ x \in V \, : \, gx=x \}.$$ Then 
the number 
$$ e(V,G) := \frac{1}{|G|} \sum_{ \begin{array}{c} {\scriptstyle (g,h) 
\in G\times G} \\
{\scriptstyle  gh=hg }\end{array} } e(V^g \cap V^h) $$
is called the {\bf physicists' orbifold Euler number} of $V$.
}
\label{orb-ph}
\end{dfn}

Our main result in this paper is the following: 

\begin{theo}
Let $V$ be as in \ref{orb-ph},  $X:= V/G$ the geometric quotient, 
$\Delta_1, \ldots, \Delta_k \subset V/G$ the set of all  
irreducible components of codimension $1$ 
in the ramification locus of the Galois covering 
$\phi\,:\, V \rightarrow X$. We denote by  $\nu_i$ the order of 
a  cyclic inertia subgroup $G_i \subset G$ corresponding to 
$\Delta_i$   and 
set 
\[ \Delta_X  : = \sum_{i =1}^k \left(\frac{\nu_i-1}{\nu_i}\right) \Delta_i. \]
Then the pair $(X, \Delta_X)$ is Kawamata log terminal and 
the following equality holds
\[  e_{\rm st}(X, \Delta_X) = e(V, G). \]
\label{equals}
\end{theo} 

As an corollary of \ref{equals}, we obtain the following 
statement  conjectured by  Miles Reid in \cite{R}: 

\begin{theo} 
Let $G \subset {\rm SL}(n, \C)$ be a finite subgroup acting on $V: = \C^n$.
Assume that there exists a crepant desingularization of $X:= V/G$, i.e.,
a smooth variety $Y$ together with a projective birational morphism
$\rho\, : \, Y \rightarrow X$ such that the canonical class  $K_Y$ 
is trivial. Then the Euler number of $Y$ equals the number of conjugacy 
classes in $G$. 
\label{equals2}
\end{theo}

The paper is organized as follows. In Section 2  we review a construction of 
a non-Archimedian measure on the space of arcs $J_{\infty}(X)$ of 
a smooth algebraic variety $X$ over $\C$. 
This measure associate to 
a measurable subset $C \subset J_{\infty}(X)$ an element $Vol_X(X)$ 
of a  $2$-dimensional notherian ring $\widehat{A}_1$ which is complete 
with respect to a non-Archimedian topology defined by powers of a principal 
ideal $( \theta) \subset \widehat{A}_1$. In Section 3 we define   
expenentially integrable measurable functions and their  
exponential non-Archimedian integrals. 
Our main interest are measurable functions $F_D$ associated with 
$\Q$-divisors $D \in Div(X) \otimes \Q$. We prove \ref{indep} using a 
transformation formula for the exponential integral under a birational 
proper morphism. 

In Section 4 we consider Kawamata log terminal 
pairs $(X, \Delta_X)$,  where $X$ is a toric variety and $\Delta_X$ 
is a torus invariant $\Q$-divisor. We give an explicit formula 
for $e_{\rm st}(X, \Delta_X)$  using a $\Sigma$-piecewise linear function 
$\varphi_{K,\Delta}$ corresponding  to the torus invariant  
$\Q$-Cartier divisor 
$K_X + \Delta_X$. In Section 5 we investigate quotients of smooth algebraic 
varieties $V$ modulo regular actions of finite groups $G$. We define    
canonical sequences of blow ups of smooth $G$-invariant subvarieties 
in $V$ which allow us to construct  
in a canonical way  a smooth $G$-variety $V'$ 
such that stabilizers 
of all points in $V'$ are abelian. This construction is used in Section 6 
where we prove  our  main theorem \ref{equals}. In Section 7 we apply our 
results to a cohomological McKay correspondence in arbitrary dimension 
(this extends our $p$-adic ideas from \cite{B}). 

We note that Sections 2 and 3 are strongy influenced by the idea of ``motivic 
integral'' proposed by Kontsevich \cite{K}. Its different versions 
are containend in the papers of  Denef and Loeser \cite{D2,DL,DL1,DL2}.  
The case of divisors on surfaces was considered  by Veys in \cite{V1,V2}. 

\bigskip
\bigskip

\noindent
{\bf Acknowledgements:}$\;\;$  It is my  pleasure to thank Professors
Yujiro  Kawamata, Maxim  Kontsevich, Shigefumi Mori, and Miles Reid   
for useful  discussions.

\section{Non-Archimedian measure on spaces of arcs}

Recall definitions of jets and  spaces of arcs 
(see  \cite{G-G}, Part A).

\begin{dfn} 
{\rm Let $X$ be a smooth $n$-dimensional complex manifold, $x \in X$ an 
arbitrary point. A {\bf germ of a holomorphic curve at $x$} is  
a germ of a holomorphic map $\gamma$ of a small ball $\{ |z|< \varepsilon \} \subset \C$ 
to $X$ such that $\gamma(0) =x$.  

Let $l$ be a nonnegative integer. 
Two germs $\gamma_1, \gamma_2$ of holomorphic curves at $x$ 
are called {\bf $l$-equivalent} 
if the derivatives of $\gamma_1$ and $\gamma_2$ at  $0$ coincide up 
to order $l$. The set of 
$l$-equivalent germs of holomorphic curves is  denoted by 
$J_l(X,x)$ and called the {\bf jet space  of order $l$ at $x$}. 
The union 
$$J_l(X) = \bigcup_{x \in X} J_l(X,x)$$  
is a complex manifold of dimension 
$(l+1)n$ which is a holomorphic affine  bundle over $X$.  
The complex manifold 
$J_l(X)$ is  called the {\bf jet space  of order $l$ of  $X$}. } 
\end{dfn}

\begin{dfn} 
{\rm Consider canonical mappings  
$j_l\, : \, 
J_{l+1}(X) \to J_l(X)$ $(l \geq 0)$ whose fibers 
are isomorphic to affine spaces $\C^n$. 
We denote by 
$J_{\infty}(X)$ the projective limit of $J_l(X)$ and by 
$\pi_l$ the canonical projection $J_{\infty}(X) \to J_l(X)$. 
The space  $J_{\infty}(X)$  is called the {\bf space of arcs of $X$}.    
} 
\end{dfn} 

\begin{rem}  
{\rm Let  $R$ be the formal power series ring $\C[[t]]$ considered as 
the inverse limit of 
finite dimensional $\C$-algebras $R_l: = \C[t]/(t^{l+1})$.
If  $X$ is  $n$-dimensional  smooth quasi-projective algebraic 
variety over $\C$, then the set of points in 
$J_{\infty}(X)$ (resp. $J_l(X)$) coincides with the set of   
$R$-valued (resp. $R_l$-valued) points of $X$. } 
\end{rem} 

From now on  we shall consider only the spaces $J_{\infty}(X)$, 
where $X$ is a smooth algebraic variety. In this 
case, $J_l(X)$ is a smooth  algebraic variety 
for all $l \geq 0$.

\begin{dfn} 
{\rm A set $C \subset J_{\infty}(X)$ is called {\bf cylinder set} if 
there exists a positive integer $l$ such that 
$C = \pi^{-1}_l(B_l(C))$ for some constructible subset    
$B_l(C) \subset J_l(X)$. Such a constructible subset  $B_l(C)$ 
will be called the 
$l$-{\bf base of} $C$. By definition, the empty 
set $\empty \subset  J_{\infty}(X)$
is a cylinder set  and its $l$-base in $J_l(X)$
is assumed to be empty for all $l \geq 0$.} 
\end{dfn} 

\begin{rem} 
{\rm  Let $C  \subset J_{\infty}(X)$ be a cylinder set with 
an $l$-base $B_l(X)$. 

(i) It is clear that $B_{l+1}(C):= 
j^{-1}_l(B_l(C)) \subset J_{l+1}(X)$ is the $(l+1)$-base of $C$ and 
$B_{l+1}(X)$ is a Zariski locally trivial affine bundle over 
$B_l(C)$  whose fibers are 
isomorphic to $\C^n$. 

(ii)  Using (i), it is a standard exercise to show that finite unions, 
intersections  and complements  of 
cylinder sets are  again cylinder sets. 
\label{cyl}
} 
\end{rem}

The following property of cylinder sets will be important: 

\begin{theo} 
Assume that a cylinder set $C \subset J_{\infty}(X)$ is contained 
in a countable union $\bigcup_{i =1}^{\infty} C_i$ of 
cylinder sets $C_i$. Then there exists a positive integer  $m$ such that 
$C \subset \bigcup_{i =1}^{m} C_i$. 
\label{cover}
\end{theo} 

\noindent
{\em Proof.} The proof of theorem \ref{cover} 
is based on a classical property of 
constructible sets (see \cite{G-D}, Cor. 7.2.6). For  details see
Theorem 6.6 in \cite{B1}. Another version of the same statement  
is contained in  \cite{DL1} (see Lemma 2.4). 
\hfill $\Box$

\begin{dfn} 
{\rm Let $\Z[\tau^{\pm 1}]$ be the Laurent polynomial 
ring in variable $\tau$ with coefficients in $\Z$,   
${A}$ the group algebra 
of $(\Q, +)$ with coefficients in  $\Z[\tau^{\pm 1}]$. 
We denote by $\theta^s \in {A}$ the image of 
$s \in \Q$ under the natural  homomorphism 
$(\Q, +) \to  ({A}^*, \cdot)$, where  ${A}^*$ 
is the multiplicative group of invertible elements in ${A}$ 
(the element $\theta \in {A}$ is transcendental over 
$\Z[\tau^{\pm 1}]$). For this reason, we write  
$${A}:=  \Z[ \tau^{\pm 1}][\theta^{\Q}] $$ 
and  identify ${A}$  with the  direct limit of the subrings
${A}_N :=  \Z[ \tau^{\pm 1}][\theta^{\frac{1}{N}\Z}] \subset 
{A}$, where $N$ runs over all positive integers.
} 
\end{dfn}  

\begin{dfn} 
{\rm We consider a topology on ${A}$ defined by 
the {\bf non-Archimedian norm}
\[ \| \cdot \| \; : \; {A} \to \R_{\geq 0}  \]
which is uniquely characterised by the properties: 

(i)  $\| ab \|  = \| a\| \cdot \| b\|$, 
$\forall a,b \in  {A}$; 

(ii)  $\| a + b \|  =  \max \{ 
 \| a\|, \| b \| \}$, 
$\forall a,b \in  {A}$ if  $\| a\| \neq  \| b \|$; 

(ii)  $\| a \|  = 1$, $\forall a \in \Z[ \tau^{\pm 1}] 
\setminus \{0 \}$; 

(iii) $\| \theta^s \|  = e^{-s}$ if $s \in \Q$. 

\noindent
The {\bf completion} of ${A}$ (resp. of ${A}_N$) with respect 
to the norm $\| \cdot \|$ will be denoted  
by $\widehat{A}$ (resp. by  $\widehat{A}_N$).  We set 
\[ \widehat{A}_{\infty} := \bigcup_{N \in \N} 
\widehat{A}_N \subset \widehat{A}. \]}
\end{dfn}  

\begin{rem}
{\rm  The noetherian ring 
$\widehat{A}_N$ consists of Laurent power series in variable 
$\theta^{1/N}$ with 
coefficients in $\Z[ \tau^{\pm 1}]$. The ring  $\widehat{A}$ 
consists consists of formal infinite sums
\[ \sum_{i=1}^{\infty} a_i \theta^{s_i},\;\;a_i \in \Z[ \tau^{\pm 1}], \]
where $s_1 <  s_2 < \cdots $ is an ascending sequence of rational numbers 
having the property $\lim_{i \to \infty} s_i = +\infty$.} 
\end{rem} 

\begin{dfn} 
{\rm Let $W$ be an arbitrary algebraic  variety. 
Using a natural mixed Hodge structure in cohomology 
groups   $H^i_c(W, \C)$, $( 0 \leq i \leq 2d)$, we define  
the number $h^{p,q}\left(H^i_c(W, \C)\right)$ to be the dimension 
of the $(p,q)$-type Hodge component in $H^i_c(W, \C)$. 
We set 
\[ e^{p,q}(W): = 
\sum_{i \geq 0} (-1)^i h^{p,q}\left(H^i_c(W, \C)\right) \]
and call 
\[ E(W; u,v) := \sum_{p,q} e^{p,q}(W) u^p v^q,  \]
the {\bf $E$-polynomial} of $W$.  
By the {\bf usual Euler number} 
of $W$ we always mean 
$ e(W) := E(W;1,1)$. } 
\end{dfn}

\begin{rem} 
{\rm For our purpose, it will be very important that $E$-polynomials 
have properties which are very similar to the ones of  usual Euler numbers: 

(i) if  $W = W_1 \cup \cdots \cup W_k$ is a disjoint union of 
Zariski locally closed subsets  $W_1, \ldots, W_k$,  then 
\[ E(W; u,v) = \sum_{i=1}^k E(W_i; u,v); \]  

(ii) if $W = W_1 \times W_2$ is a product of two algebraic varieties 
$W_1$ and $W_2$, then 
\[ E(W; u,v) = E(W_1; u,v) \cdot E(W_2; u,v); \]

(iii) if $W$ admits a fibering over $Z$  which is locally
trivial in Zariski topology such that each fiber of the morphism 
$f\,: \, W \to Z$ is isomorphic to the affine space ${\C}^n$, then 
\[ E(W; u,v) =  E({\C}^n; u,v) \cdot E(Z; u,v)= (uv)^n E(Z; u,v). \]
}
\label{e-poly} 
\end{rem} 
 
\begin{dfn} 
{\rm Let  $V \subset W$ is a constructible subset in a 
complex algebraic variety $V$. We write $V$ as a union  
\[ V = W_1 \cup \cdots \cup W_k \]
of pairwise nonintersecting Zariski locally closed subsets 
$W_1, \ldots, W_k$. Then the  {\bf $E$-polynomial of} $V$ is defined as 
follows:
\[ E(V; u,v) := \sum_{i=1}^k E(W_i; u,v). \]} 
\label{const-e}
\end{dfn} 

\begin{rem} 
{\rm Using \ref{e-poly}(i), it is easy to check that the above definition
does not depend on the choice of the decomposion of $V$ into a finite 
union of pairwise nonintersecting  Zariski
locally closed subsets.} 
\end{rem}

Now we define a {\bf non-Archimedian cylinder set measure} 
on $J_{\infty}(X)$.

\begin{dfn} 
{\rm  $C \subset J_{\infty}(X)$ be a  cylinder set. 
We define the {\bf non-Archimedian volume} $Vol_X(C) \in 
{A}_1$  of $C$ by the 
following formula: 
\[ Vol_X(C) := 
E(B_l(C); \tau \theta^{-1}, \tau^{-1} \theta^{-1}) \theta^{(l+1)n} \in 
{A}_1,  \]
where  $C = \pi^{-1}_l(B_l(C))$ and 
$E(B_l(C); u,v)$ is the $E$-polynomial of the $l$-base 
$B_l(C) \subset J_l(X)$. If $C= \emptyset$, we set $Vol_X(C) :=0$.   
}
\end{dfn}

\begin{rem} 
{\rm Using \ref{cyl}(i) and \ref{e-poly}, one immediately obtains that 
$Vol_X(C)$ does not depend on the choice of 
an $l$-base $B_l(C)$ and 
$$\|Vol_X(C)\| = e^{2dim\, B_l(C)  - 2(l+1)n}.$$  
In particular, one has the following properties

{\rm (i)} If $C_1$  and $C_2$ are two cylinder sets such that $C_1 
\subset C_2$, then 
\[ \| Vol_X(C_1) \| \leq  \| Vol_X(C_2) \|. \]

{\rm (ii)} If $C_1, \ldots, C_k$ are cylinder sets, then 
\[ \| Vol_X(C_1 \cup \cdots \cup C_k) \| =   \max_{i =1}^k 
\| Vol_X(C_i) \|. \]

{\rm (iii)} if a cylinder set 
$C$ is a finite disjoint union of cylinder sets $C_1, \ldots, 
C_k$, then 
\[ Vol_X(C)= Vol_X(C_1) + \cdots + Vol_X(C_k). \]
\label{prop12} } 
\end{rem}

\begin{dfn} 
{\rm  We say that a subset $C \subset J_{\infty}(X)$ 
is  {\bf measurable} if  for any positive real number $\varepsilon$
there exists   a sequence of cylinder sets  $C_0(\varepsilon), 
C_1(\varepsilon), C_2(\varepsilon), \cdots $ such that 
$$ \left( C \cup  C_0(\varepsilon) \right) \setminus 
\left( C \cap  C_0(\varepsilon) \right) 
\subset \bigcup_{i \geq 1} C_i(\varepsilon) $$ 
and $\| Vol_X(C_i({\varepsilon}))\| < \varepsilon$ for all $i \geq 1$.
If $C$ is measurable, then the element 
\[ Vol_X(C) := \lim_{\varepsilon \to 0} C_0(\varepsilon) 
\in \widehat{A}_1 \]
will be called the {\bf non-Archimedian volume} of $C$. 
}
\end{dfn}

\begin{theo} 
If $C \subset J_{\infty}(X)$ is measurable, then 
$\lim_{\varepsilon \to 0} 
C_0(\varepsilon)$
exists and  does not depend on the choice of sequences   
$C_0(\varepsilon), C_1(\varepsilon),  C_2(\varepsilon), \cdots $.   
\label{lim}
\end{theo} 

\noindent
{\em Proof. } The property \ref{cover} plays a crucial 
role in the  proof of this theorem. For details see 
\cite{B1}, Theorem 6.18. 
\hfill $\Box$  

The proof of the following statement is a standard exercise: 

\begin{prop} 
Measurable sets possess the following properties: 

{\rm (i)}  Finite unions, finite intersections  of 
measurable sets are measurable. 

{\rm (ii)} If $C$ is a disjoint union of nonintersecting measurable 
sets $C_1, \ldots, C_m$, then 
\[ Vol_X(C) = Vol_X(C_1) + \cdots + Vol_X(C_m). \] 

{\rm (iii)} If $C$ is measurable, then the complement 
$\overline{C}:= J_{\infty}(X) \setminus C$ is measurable.

{\rm (iv)}   If $C_1, C_2, \ldots, C_m, \ldots$ is an infinite sequence 
of nonintersecting measurable sets having the property 
\[ \lim_{i \to \infty} \| Vol_X(C_i) \| =0, \]
then 
\[ C = \bigcup_{i =1}^{\infty} C_i \]
is measurable and 
\[ Vol_X(C) = \sum_{i =1}^{\infty} Vol_X(C_i). \]   
\label{bool}
\end{prop} 

The next example shows that our  non-Archimedian 
measure does not have all  properties of  the standard Lebesgue 
measure: 

\begin{exam} 
{\rm Let $C \subset R= \C[[t]]$ be the set consisting of 
all power series $\sum_{i \geq 0} a_i t^i$ such that $a_i \neq 0$ 
for all $i \geq 0$.  
For any $k \in \Z_{\geq 0}$, we define $C_k \subset R$ to 
be the set consisting of 
all power series $\sum_{i \geq 0} a_i t^i$ such that $a_i \neq 0$ 
for all $0 \leq i \leq k$. We identify $R$ with $J_{\infty}({\C})$. 
Then every $C_k \subset  J_{\infty}({\C})$ is a cylinder set 
and  $Vol_{{\C}}(C_k)= (1 -\theta^2)^{k+1}$. Moreover, we have 
\[ C_0 \supset C_1 \supset C_2 \supset \cdots\; , \;\; \mbox{\rm and} \;\; 
C = \bigcap_{k\geq 0} C_k. \]
However, the sequence
\[   Vol_{{\C}}(C_0), \; Vol_{{\C}}(C_1), \;  
Vol_{{\C}}(C_2), \; \ldots \]
does not converge in $\widehat{A}_1$.} 
\end{exam}

\begin{dfn} 
{\rm We shall say that a subset $C \subset J_{\infty}(X)$ 
has  {\bf measure zero} 
if for any positive real number $\varepsilon$ there exists a sequence 
of cylinder sets $C_1(\varepsilon), C_2(\varepsilon), \cdots $ such that 
 $C \subset \bigcup_{i \geq 1} C_i(\varepsilon)$ and  
$\| Vol_X(C_i({\varepsilon}))\| < \varepsilon$ for all $i \geq 1$. } 
\end{dfn}

\begin{dfn} 
{\rm Let $Z \subset X$ be a Zariski closed subvariety.   
For any point $x \in Z$, we denote by ${\cal O}_{X,x}$ the ring 
of germs of holomorphic functions at  $x$. 
Let $I_{Z,x} \subset {\cal O}_{X,x}$ be the ideal of germs of holomorphic 
functions vanishing on $Z$. We set 
\[ J_{l}(Z,x) := \{ y \in J_{l}(X,x)\, : \, g(y) =0 \;\; 
\forall \; g \in  I_{Z,x} \}, \;\; l \geq 1,  \]
\[ J_{\infty}(Z,x) := \{ y \in J_{\infty}(X,x)\, : \, g(y) =0 \;\; 
\forall \; g \in  I_{Z,x} \} \]
and 
\[   J_{\infty}(Z) := \bigcup_{x \in Z}  J_{\infty}(Z,x). \]
The space 
$J_{\infty}(Z) \subset J_{\infty}(X)$ will be called {\bf space of arcs 
with values in $Z$}. 
} 
\end{dfn} 

\begin{prop} 
Let  $Z$ be an arbitrary Zariski closed subset in 
a smooth irredicible algebraic variety $X$.
Then $J_{\infty}(X,Z) \subset J_{\infty}(X)$ is measurable. Moreover, 
one has  
\[ Vol_X(J_{\infty}(Z)) = \left\{ \begin{array}{ll} 0
\; & \; 
\mbox{\rm if $Z \neq  X$} \\ 
Vol_X(J_{\infty}(X))\; &\; 
\mbox{\rm if $Z = X$.} 
\end{array} \right. \]  
\label{subv}
\end{prop} 

\noindent
{\em Proof.} If $Z \neq X$, then the set $J_{\infty}(Z)$ can be obtained 
as an intersection of cylinder sets $C_k$ such that $\| Vol_X(C_k) \| \leq 
e^{-2k}$ (see Theorem 6.22 in \cite{B1} and  3.2.2 in \cite{DL1}). 
\hfill $\Box$ 

\section{Non-Archimedian integrals} 

\begin{dfn}
{\rm  By a {\bf measurable function $F$ on} $ J_{\infty}(X)$ we mean  
a function  $F\, : \, M  \to  
\Q$, where $M \subset  J_{\infty}(X)$ is  a subset such that 
$ J_{\infty}(X) \setminus M$ has measure zero and   
 $F^{-1}(s)$ is measurable 
for all  $s \in \Q$. Two measurable functions  $F_i\, : \, M_i  \to  
\Q$ $(i =1,2)$   
on  $ J_{\infty}(X)$  are called {\bf equal} 
if $F_1(\gamma) = F_2(\gamma)$ for all 
$\gamma \in M_1 \cap M_2$.
} 
\end{dfn} 

\begin{dfn} 
{\rm A measurable function  $F\, : \,  M  
\to \Q$ is called 
{\bf exponentially integrable}  if the series 
\[ \sum_{s \in \Q}  \| Vol_X(F^{-1}(s)) \| e^{-2s} \]
converges.  
 If $F$ is exponentially integrable, then 
the sum 
\[  \int_{ J_{\infty}(X)}  e^{-F} := 
  \sum_{s \in \Q}  Vol_X(F^{-1}(s)) \theta^{2s} \in \widehat{A} \]
will be called the {\bf exponential 
integral of ${F}$ over $ J_{\infty}(X) $}. }
\end{dfn}

\begin{dfn} 
{\rm Let $D \subset Div(X)$ be a  subvariety of codimension $1$, 
$x \in D$ a point, and $g \in {\cal O}_{X,x}$ the local equation for 
$D$ at $x$. We set  $M(D) :=  J_{\infty}(X) \setminus J_{\infty}(Supp\, D)$.  
For any $\gamma  \in M(D)$, we denote by 
$\langle D, \gamma \rangle_x$ the order of the holomorphic function 
$g(\gamma(t))$ at $t =0$. The number 
$\langle D, \gamma \rangle_x$ will be called the {\bf intersection number} 
of $D$ and $\gamma$ at $x \in X$.  We define the  function 
$F_D\; : \; M(D)  \to \Z$
as follows: 
\[ F_D(\gamma) = \left\{ \begin{array}{ll} 0
\; & \; 
\mbox{\rm if $\pi_0(\gamma) = x \not\in D$} \\ 
\langle D, \gamma \rangle_x \; &\; 
\mbox{\rm if $\pi_0(\gamma) \in D$} 
\end{array} \right. \] 
}
\end{dfn}

\begin{rem}
{\rm Using the property $\langle D' + D'', \gamma \rangle_x = 
\langle D', \gamma \rangle_x + \langle D'', \gamma \rangle_x$, we extend 
the definition of $F_D$ to an arbitrary $\Q$-Cartier divisor $D$:  
if $D = \sum_{i =1}^m a_iD_i \in Div(X) \otimes \Q$ is a 
$\Q$-linear combination of irreducible subvarieties $D_1, \ldots, D_m$, then  
we set 
\[ F_D := \sum_{i =1}^m  a_i F_{D_i}.  \]
It is easy to show that measurable functions form a $\Q$-vector space and 
$D \subset Div(X) \otimes \Q$ can be identified with its  $\Q$-subspace, 
since   $F_D \,: \, M(D)  \to \Q$
is mesurable for all $D \subset Div(X) \otimes \Q$     
}
\end{rem} 

The following theorem describes a transformation law for exponential 
integrals  under proper birational morphisms: 

\begin{theo} 
Let $\rho\, : \, Y \to X$ be a proper birational  morphism of 
smooth complex algebraic varieties, $D = \sum_{i =1}^r d_i D_i \in Div(Y)$ 
the Cartier divisor defined by the equality 
\[ K_{Y} = \rho^* K_X +  \sum_{i =1}^r d_i D_i.\]
Denote by $\rho_{\infty} \,: \,  J_{\infty}(Y) \to J_{\infty}(X)$ the 
mapping of spaces of arcs induced by $\rho$. 
Then  a  measurable function  $F$ 
is exponentially
integrable if an only if $F \circ \rho_{\infty}  + F_{D}$ is exponentially  
integrable. Moreover, if the latter holds, then 
\[ \int_{J_{\infty}(X)} e^{-F} = \int_{J_{\infty}(Y)} 
e^{- F  \circ \rho_{\infty} - F_{D}}. \]
\label{morph}
\end{theo} 

\noindent
{\em Proof.}  The proof of theorem \ref{morph} is based on the equality
$Vol_Y(C) = Vol_X(\rho_{\infty}(C)) \theta^{2a}$, where $C$ is a cylinder
set in $J_{\infty}(Y)$ such that $F_D(\gamma)= a$ for all $\gamma \in C$ 
(see for details  Theorem 6.27 in  \cite{B1} and Lemma 3.3  in \cite{DL1}). 
\hfill $\Box$

\begin{theo} 
Let  $D:= a_1D_1 + \cdots + a_m D_m \in Div(X) \otimes \Q$ be 
a $\Q$-divisor. Assume $Supp\, D$ 
is a normal crossing divisor. 
Then $F_D$ is exponentially integrable if and only 
$a_i > -1$ for all $i \in \{1, \ldots, m \}$. Moreover, if the 
latter holds, then 
\[  \int_{J_{\infty}(X)} e^{-F_D} = \sum_{J \subset I}  
E(D^{\circ}_J; \tau \theta^{-1},\tau^{-1} \theta^{-1} ) 
(\theta^{-2} -1)^{|J|}\prod_{j \in J} \frac{1}{1 - \theta^{2(1 + a_j)}}
   \]
\label{int-for}
\end{theo}

\noindent
{\em Proof.} The set $M(D) \subset J_{\infty}(X)$ 
splits into a countable union of pairwise nonintersecting cylinder sets 
whose non-Archimedian volume can be computed via $E$-polynomials 
of the starta $D^{\circ}_J$ (see for details Theorem 6.28 in  \cite{B1}  
and Theorem 5.1 in \cite{DL1}). 

\hfill $\Box$

\begin{dfn} 
{\rm Let  $(X, \Delta_X)$ be a Kawamata log terminal pair. Consider 
a log resolution $\rho_1\, : \, Y \to X$ and write 
\[ K_{Y} = \rho^* (K_X + \Delta_X) + \sum_{i=1}^m a(D_i, \Delta_X) D_i.  \] 
Using the notations from \ref{ph}, we define  
\[ E_{\rm st}(X, \Delta_X) := \sum_{J \subset I} E(D_J^{\circ}; u,v) 
\prod_{j \in J} 
\frac{uv-1}{(uv)^{a_l(D_j, \Delta_X)} -1}. \]
The function $E_{\rm st}(X, \Delta_X; u,v)$ will be called 
{\bf stringy $E$-function of $(X,  \Delta_X)$}. 
} 
\end{dfn}

\begin{theo} 
Let $(X, \Delta_X)$ be a Kawamata log terminal pair. Then 
the  $E$-function of $(X, \Delta)$ does not depend on the choice of 
a log resolution. 
\label{e-st1}
\end{theo} 

\noindent
{\em Proof.} Let $\rho_1\, : \, Y_1 \to X$ and 
$\rho_2\, : \, Y_2 \to X$ be two log 
resolutions of singularities such that 
\[ K_{Y_1} = \rho_1^* (K_X + \Delta_X) + D_1, \;\; 
K_{Y_2} = \rho_2^* (K_X + \Delta_X) + D_2 \]
where 
\[ D_1 = \sum_{i=1}^{r_1} a(D_i', \Delta_X)D_i' 
\; \;\; \mbox{\rm and} \;\;  
D_2 =  \sum_{i=1}^{r_2} a(D_i'', \Delta_X)D_i'' \]
and all discrepancies $a(D_i', \Delta_X)$,    
$a(D_i'', \Delta_X)$ are $> -1$. 
Choosing a resolution of singularities  $\rho_0\, : \, Y_0\to X$ which 
dominates both resolutions $\rho_1$ and $\rho_2$, we obtain 
two morphisms $\alpha_1\, : \, Y_0\to Y_1$ and 
$\alpha_2\, : \, Y_0 \to Y_2$ such that 
$\rho_0 = \rho_1 \circ \alpha_1 =  \rho_2 \circ \alpha_2$.  
We set  $F := F_{D_0}$, where $D_0 =  K_{Y_0} -  \rho_0^* (K_X + \Delta_X)$. 
Since 
\[ K_{Y_0} -  \rho_0^* (K_X  + \Delta_X)=  
( K_{Y_0} - \alpha_i^* K_{Y_i})  + \alpha_i^*D_i, \;\; (i =1,2),  \]
we obtain 
\[  \int_{J_{\infty}(Y_1)} e^{-F_{D_1}} =  
\int_{J_{\infty}(Y_0)} e^{-F_{D_0}} =  
\int_{J_{\infty}(Y_2)} e^{-F_{D_2} }\;\;\;
\mbox{\rm (see \ref{morph})}. \]
It follows from  \ref{int-for} that  
\[ \int_{J_{\infty}(Y_i)} e^{-F_{D_i}} =
 E_{\rm st}(X,\Delta_X; \tau \theta^{-1},  
\tau^{-1} \theta^{-1}), \;\; i \in \{0,1,2\}.\]
Making the  substitutions $u = \tau \theta^{-1},  
v=  \tau^{-1} \theta^{-1}$, 
we obtain that the definition of the stringy $E$-function 
$E_{\rm st}(X, \Delta_X; u,v)$  does not depend 
on the choice of log resolutions $\rho_1$ and $\rho_2$. 
\hfill $\Box$  

\bigskip

\noindent
{\bf  Proof of Theorem \ref{indep}:} 
The statement immediately follows from \ref{e-st1} using the equality  
\[ e_{\rm st}(X, \Delta_X) = \lim_{u,v \to 1} E_{\rm st}(X, \Delta_X; u,v). \]

\hfill $\Box$

\section{Log pairs on toric varieties}

Let $X$ be a normal toric variety of dimension $n$ 
associated with a rational polyhedral 
fan $\Sigma \subset N_{\R} = N \otimes {\R} $, where  
${N}$ is a free abelian group of rank $n$. Denote by $X({\sigma})$ 
the torus orbit in $X$ corresponding to a cone $\sigma \in \Sigma$ 
($codim_X X_{\sigma} = dim\, \sigma$). 
Let  $\overline{X}({\sigma})$ be the Zariski  closure of  $X({\sigma})$.
Then the torus invariant $\Q$-divisors are $\Q$-linear combinations 
of the closed strata $\overline{X}(\sigma_1^{(1)}), \ldots, 
\overline{X}(\sigma_k^{(1)})$, where $\Sigma^{(1)} := \{ \sigma_1^{(1)}, 
\ldots, \sigma_k^{(1)} \}$ is the set of all $1$-dimensional cones in 
$\Sigma$. We denote by $e_1, \ldots, e_k$  the primitive lattice generators
of the cones $\sigma_1^{(1)}, \ldots, \sigma_k^{(1)}$ and set  
$\Delta_i := \overline{X}(\sigma_i^{(1)})$ $i \in \{1, \ldots, k \}$. 
 
\begin{dfn} 
{\rm   Let $\varphi_{K,\Delta}\, : \, N_{\R} \rightarrow 
{\R}_{\geq 0}$  be a continious function satisfying
the conditions 

(i)  $\varphi_{K,\Delta}(N) \subset  \Q$; 

(ii) $\varphi_{K,\Delta}$ is linear on each cone $\sigma \in \Sigma$; 

(iii) $\varphi_{K,\Delta} (p) >0 $  for all $p \in N \setminus \{ 0\}$. 

\noindent
Then we define a {\bf $\Q$-divisor $\Delta_X \in Z_{n-1}(X)$ associated with} 
$\varphi_{K,\Delta}$ 
as follows: 
\[ \Delta_X : = \sum_{i =1}^k \left( 1- \varphi_{K,\Delta} (e_i)  
\right) \Delta_i. \]
} 
\label{tor-l}
\end{dfn} 

\begin{rem} 
{\rm It is well-known that the canonical class $K_X$ of a toric variety 
$X$ is equal to $- (\Delta_1 + \cdots + \Delta_k)$. The above definition 
of $\Delta_X$ implies that $K_X + \Delta_X$ is a $\Q$-Cartier divisor
on $X$ corresponding to the 
$\Sigma$-piecewise linear function $-\varphi_{K,\Delta}$. } 
\end{rem} 

The following statement is well-known in toric geometry (see e.g. 
\cite{KMM} \S 5-2): 

\begin{prop} 
Let $\rho\: : \; X' \to X$ be a toric desingularization of $X$, which 
is defined by a subdivision $\Sigma'$ of the fan $\Sigma$. 
Denote by  $\{ D_1, \ldots, D_{m} \}$ the set of all irreducible torus 
invariant strata on $Y$ corresponding to  primitive lattice  
generators $ e_1', \ldots, e_{m}'$  of $1$-dimensional cones  
$\sigma' \in \Sigma'$.  
Then $\sum_{i =1}^{m} D_i$ is a normal crossing divisor and one has 
\[ K_{X'} = \rho^* (K_X + \Delta_X)  + \sum_{i=1}^{m} 
a(D_i, \Delta_X)D_i,    \]
where  $a(D_i, \Delta_X)  = 
\varphi_{K,\Delta}(e_i') -1$ $\forall i \in \{1, \ldots, m\}$. 
\label{tor-klt}
\end{prop} 

\begin{coro} 
Let $\varphi_{K,\Delta}$ be a $\Sigma$-piecewise linear function as in 
\ref{tor-l}. Then the pair $(X, \Delta_X)$ is Kawamata log termial. 
\label{tor-klt1}
\end{coro} 
      
Denote by   $\sigma^{\circ}$ the relative interior of $\sigma$ 
(we put $\sigma^{\circ} = 0$,
 if $\sigma = 0$).  We give the following  explicit formula for the 
function $E_{\rm st}(X, \Delta_X; u, v)$:  

\begin{theo} 
\[ E_{\rm st}(X, \Delta_X; u,v) = (uv -1)^n \sum_{ \sigma \in \Sigma} 
\sum_{p \in \sigma^{\circ} \cap N} (uv)^{-\varphi_{K,\Delta}(p)} =
(uv -1)^n  \sum_{p \in N} (uv)^{-\varphi_{K,\Delta}(p)} .  \] 
\label{st-tor}
\end{theo} 

\noindent
{\em Proof.} Let $T \subset X$ be an algebraic torus acting on $X$,  
$\partial X := X \setminus T$ its complement. 
Choose an isomorphism $N \cong \Z^n$ and write $p= (p_1, \ldots, p_n) \in 
\Z^n$. Denote by $K:= \C((t))$ the field of Laurent power series and 
define a cylinder subset $C_p \subset  J_{\infty}(X)$ as follows:
\[ C_p := \{ (x_1(t), \ldots, x_n(t)) \in K^n \;  : \; 
Ord_{t=0} x_i(t) = p_i, \; 1 \leq i \leq n \}. \] 
Consider the 
subset $M(\partial X) \subset J_{\infty}(X)$ consisting of all arcs which are 
not contained in  $J_{\infty}(\partial X)$. Then   $M(\partial X)$ 
splits into a disjoint union 
\[ M(\partial X) = \bigcup_{p \in N} C_p. \] 
Let $\rho\: : \; X' \to X$ be a toric desingularization of $X$, and 
\[ K_{X'} = \rho^* (K_X + \Delta_X)  + \sum_{i=1}^{m} a(D_i, \Delta_X)D_i. \]
By definition, we have 
\[  E_{\rm st}(X, \Delta_X; \tau \theta^{-1}, \tau^{-1} \theta^{-1}) 
= \int_{J_{\infty}(X')} e^{-F_D}, \]
where 
\[  D=  \sum_{i=1}^{m} 
a(D_i, \Delta_X)D_i. \]
Now we notice that $F_D$ is constant on each cylinder set $C_p$ $(p \in N)$ 
and 
\[ Vol(C_p) \theta^{2F(C_p)} = (\theta^{-2} -1)^n 
\theta^{2\varphi_{K,\Delta}(p)}. \]
Summing over $p \in N$ and making the substitution $u = \tau \theta^{-1}$, 
$v = \tau^{-1} \theta^{-1}$, we come to the required formula. \hfill $\Box$   

\begin{dfn} 
{\rm Let $X$ be an arbitary $n$-dimensional  
normal toric variety defined by a 
fan $\Sigma$, and  $X + \Delta_X$  a torus invariant 
 $\Q$-Cartier divisor corresponding to a 
$\Sigma$-piecewise linear function $\varphi_{K,\Delta}$.  
 Denote by $\Sigma^{(n)}$ the set of all $n$-dimensional 
cones in $\Sigma$.  
Let $\sigma \in \Sigma^{(n)}$ be  a cone. 
Define {\bf $\Delta$-shed of $\sigma$}
to be the pyramid
\[ {\rm shed}_{\Delta} \sigma =  \sigma \cap 
\{ y \in N \otimes \R \, : \, \varphi_{K, \Delta}(y) \leq 1 \}.  \]
Furthermore, define 
{\bf $\Delta$-shed of $\Sigma$} to be 
\[ {\rm shed}_{\Delta} \Sigma = \bigcup_{\sigma \in  \Sigma^{(n)}}  
{\rm shed}_{\Delta} \sigma. \]
}
\end{dfn}

\begin{dfn} 
{\rm Let $\sigma \in \Sigma^{(n)}$ be an arbitrary cone. 
Define $vol_{\Delta}(\sigma)$ to be the volume of ${\rm shed}_{\Delta} 
\sigma$ 
with respect to the lattice $N \subset N_{\R}$ multiplied by $n!$.  
We set 
\[  vol_{\Delta}(\Sigma) := \sum_{\sigma \in  \Sigma^{(n)}}  
vol_{\Delta}(\sigma). \]
} 
\end{dfn}

\begin{dfn} 
{\rm Let $X_0$, $X$, $X^+$ be 
$n$-dimensional normal projective toric varieties. 
Denote by $\Sigma$ (resp. by $\Sigma^+$) the fan defining 
$X$ (resp. $X^+$). Let  $(X, \Delta_{X})$ (resp. $(X^+, \Delta_{X^+})$)
be a torus invariant Kawamata 
log terminal pair defined by a  $\Sigma$-piecewise linear 
(resp.  $\Sigma^+$-piecewise linear)  function  
$\varphi_{K,\Delta}$ (resp. $\varphi_{K,\Delta}^+)$. 
Assume that we are given two equivariant projective birational toric 
morphisms $\alpha \,:\, X \to X_0$ and $\beta \, :\, X^+  \to X_0$ such that 
$- (K_X + \Delta_X)$ is $\alpha$-ample , $K_{X^+} + \Delta_{X^+}$ 
is $\beta$-ample, and both $\alpha$ and $\beta$ 
are isomorphisms in codimension $1$. Then the birational 
rational map $\psi:= {\beta}^{-1} \circ \alpha \, : \, (X, \Delta_X) 
 \dasharrow (X^+, \Delta_{X^+})$ 
is called a {\bf toric log flip} with respect to a $\Q$-Cartier divisor 
$K_X  +\Delta_X$. } 
\label{lflip1} 
\end{dfn}

\begin{prop} 
Let  $\psi \, : \, (X, \Delta_X) \dasharrow (X^+, \Delta_{X^+})$ be a toric 
log $(K_X  +\Delta_X)$-flip as above. Then 
\[  vol_{\Delta}(\Sigma) > vol_{\Delta}(\Sigma^+). \] 
\end{prop}

\noindent
{\em Proof.} Using a toric interpretation of 
ampleness via a combinatorial convexity, one obtains  
from the definition of toric log flips that 
$\varphi_{K,\Delta}(p) \leq  \varphi_{K,\Delta}^+(p)$ for all 
$p \in N$ and there exists a $n$-dimensional cone $\sigma \in 
\Sigma^{(n)}$ such that $\varphi_{K,\Delta}(p) <  
\varphi_{K,\Delta}^+(p)$ for all interior lattice points $p 
\in \sigma \cap N$. This implies the statement (cf. \cite{B1}, Prop. 4.9). 

\hfill $\Box$

\begin{prop} 
Let $X$ be an arbitary $n$-dimensional  
normal toric variety defined by a 
fan $\Sigma$, and  $X + \Delta_X$  a torus invariant 
 $\Q$-Cartier divisor corresponding to a 
$\Sigma$-piecewise linear function $\varphi_{K,\Delta}$. Then 
\[ e_{\rm st}(X, \Delta_X) =  vol_{\Delta}(\Sigma). \]
\end{prop} 

\noindent
{\em Proof.} The statement follows from the formula in \ref{st-tor} using 
the same arguments as in the proof of Prop. 4.10 in \cite{B1}. 
\hfill $\Box$

\begin{coro} 
Let $(X, \Delta) \dasharrow (X^+, \Delta_{X^+})$ be a toric log flip. 
Then 
\[ e_{\rm st}(X, \Delta_X) > e_{\rm st}(X^+, \Delta_{X^+}). \]
\label{lflip2}
\end{coro}

\section{Canonical abelianization} 

Let $G$ be a finite group, $V$ a smooth $n$-dimensional 
algebraic variety over ${\C}$  
having  a regular effective 
action of $G$. If $x \in V$ is an arbitrary 
 point,  then by $St_G(x)$ we denote the stabilizer of $x$ in $G$.
 For any element $g \in G$ we 
set $V^g := \{ x \in V\, : \, gx=x \}$.

\begin{dfn} 
{\rm Let  $D = \sum_{i=1}^m d_i D_i \in {\rm Div}(V)^G \otimes \Q$ an 
effective  
$G$-invariant $\Q$-divisor on a $G$-manifold $V$.  A pair  
$(V,D)$ will be called  {\bf $G$-normal}  if the following 
conditions are satisfied: 

(i) $Supp\, D$ is a union of normal crossing divisors $D_1, \ldots, D_m$; 

(ii) for any element 
 $g \in G$ and any irredicible component 
 $D_i$ of $D$, the divisor  $D_i$ is $St_G(x)$-invariant for 
all $x \in V^g \cap D_i$ (i.e., $h(D_i) = D_i$  $\forall \, h \in 
St_G(x)$, but the $St_G(x)$-action on $D_i$ itself 
may be nontrivial).  
} 
\end{dfn}

\begin{theo} \label{abel1}
Let $(V,D)$ be a $G$-normal pair. Then, using a 
canonically determined sequence of blow ups of $G$-invariant submanifolds, 
one obtains  a  $G$-normal pair $(V^{ab}, D^{ab})$ 
and  a projective birational $G$-morphism
$\psi\, :\, V^{ab} \rightarrow V$ 
having the properties: 

{\rm (i)} $D^{ab} = (K_{V^{ab}} -  \psi^*K_{V})  +  \psi^*D$; 

{\rm (ii)} for any point $x \in V^{ab}$  
the stabilizer $St_G(x)$ is an abelian subgroup in $G$. 
\end{theo}

\noindent 
{\em Proof.} 
Let $Z(V,G) \subset  V$ be the set of all points 
$x \in V$ such that  $St_G(x)$ is not abelian. If $Z(V,G)$ is empty, 
then we are 
done.  Assume that 
$Z(V,G) \neq \emptyset$. 
We set 
\[ s(V,G) := \max_{x \in Z(V,G)} |St_G(x)|. \]
Consider  a Zariski closed subset 
\[ Z_{\rm max}(V,G) := \{ x \in Z(V,G) \; :\;  |St_G(x)| = s(V,G) \} \subset 
Z(V,G). \]
We claim that the set  $Z_{\rm max}(V,G) \subset V$ is a smooth  $G$-invariant
subvariety of  codimension at least $2$. By definition, 
$Z_{\rm max}(V,G)$ is a  union of smooth subvarieties  
$$F(H):= \{ x \in V\; : \; gx =x\;\; \forall g \in H \},$$ 
where $H$ runs over all nonabelian 
subgroups of $G$ such that  $|H| =  s(V,G)$. 
This implies that $ Z_{\rm max}(V,G)$ is $G$-invariant.  
Since the $G$-action is effective and $dim\, F(H) =n-1$
is possible only for cyclic subgroups $H \subset G$, we obtain  
$dim\, Z_{\rm max}(V,G) \leq n-2$.  
It remains to observe that  any two subvarieties $F(H_1), F(H_2) \subset 
V$ must either coincide, or have empty 
intersection. Indeed, if $x \in F(H_1) \cap  F(H_2)$, then 
$H_1, H_2 \subset St_G(x)$. Since $|H_1|, |H_2|$ are maximal, we obtain 
$H_1 = H_2 =  St_G(x)$; i.e., $F(H_1) = F(H_2)$.

We set $V_0 := V$, $D_0 := D$ and define  $V_1$ to be the $G$-equivariant 
blow-up of $V_0$ 
with center $Z_{\rm max}(V,G)$. Denote by   
$\varphi_1 \, :\, V_1 \rightarrow V_0$ the corresponding projective birational 
$G$-morphism. It is obvious that the support of $D_1 = 
K_{V_1} - \varphi_1^*(K_V - D)$ is a normal crossing divisor. If  
$x \in V_1^g \cap E$, where $E$ is a connected component of an 
$\varphi_1$-exceptional divisor, 
then $St_G(x) \subset St_G(\varphi(x))$. Since 
$\varphi(E)$ is a connected component of a smooth subvariety 
$Z_{\rm max}(V,G)$, $ \varphi(E)$ must be  $St_G(\varphi(x))$-invariant. 
Hence, we conclude that $(V_1,D_1)$ is 
a $G$-normal pair. If $Z(V_1,G) = \emptyset$, then we are done. 
Otherwise we apply the same procedure  to the $G$-normal pair  
$(V_1, D_1)$, where $D_1 =   \phi_1^*D_0$, and  construct 
in the same way a next $G$-equivariant blow-up $\varphi_2\; : \; V_2 
\rightarrow V_1$ $\ldots$ etc.    

It remains to show that the above procedure terminates. For this purpose, it 
suffices to show that $s(V_i, G) < s(V_0, G)$ for some $i >0$. 
Assume that $s(V_0, G) =  s(V_i, G)$ for all $i >0$. Then there exist 
points $x_i \in V_i$ $(i \geq 0)$ 
such that $\varphi_{i}(x_i) = x_{i-1}$ and $St_G(x_i) = St_G(x_{i-1})$ 
$(i \geq 1)$. Let  $S(x_i)$ be the set  of those irreducible components 
of $Supp\, D_i$ which are $St_G(x_i)$-invariant and  contain  $x_i$. 
We denote by $n(x_i)$ the cardinality of $S(x_i)$ and denote by $D(x_i) 
\subset V_i$ the intersection of all divisors from $S(x_i)$. Then 
$F(St_G(x_i)) \subset D(x_i)$. If $F(St_G(x_i)) \neq D(x_i)$, then 
the point $n(x_{i+1}) = n(x_i) + 1$ (we obtain one more component from 
the $\varphi_i$-exceptional divisor over  $F(St_G(x_i))$). 
Since $n(x_i) \leq n$ for all $i \geq 0$, there exists a positive 
number $k$ such that $n(x_k) = n(x_{k+j})$ for all $j \geq 0$. 
So we obtain  $F(St_G(x_{k+j})) =  D(x_{k+j})$ for all $j \geq 0$.  
The latter means that the action of 
$St_G(x_k)$ on the tangent space to $x_k$ in $V_k$ splits into a direct 
sum of $n(x_k)$ $1$-dimensional representations and a 
$(n- n(x_k))$-dimensional trivial representation. Since the action 
of $St_G(x_k)$ is effective, the group $St_G(x_k)$ must be  abelian. 
Contradiction.  \hfill $\Box$

\begin{dfn} \label{abel3}
{\rm Let $({V}, {D})$ be a  $G$-normal pair. Then the 
$G$-normal pair  $(V^{ab},D^{ab})$ obtained in \ref{abel1} 
 will be called {\bf   canonical abelianization} of a  
$G$-normal pair $(V,D)$.}
\end{dfn}

\begin{rem}
{\rm If the stabilisator $St_G(x) \subset G$ of every  point $x \in V$ 
is already  abelian, then one can't expect that 
$G$-equivariant blow ups of smooth subvarieties  $Z \subset V$  
could  simplify singularities of 
the quotient-space $V/G$. 

Here is  the following simplest example: Let  $V := {\C}^2$ and 
$G = \langle g \rangle$ is a cyclic group of order $5$ whose 
generator $g$ acts by the diagonal matrix with the eigenvalues  
$e^{2\pi \sqrt{-1}/5}$, $e^{4\pi  \sqrt{-1}/5}$. Let $V'$ be the blow up 
of $\C^2$ at $0$. Then $V'$ has a natural covering by 
two open subsets $V_1'$ and $V_2'$ such that $V_1' \cong V_2' \cong 
\C^2$ and  the $G$-action on one of these subsets  
coincides with the original $G$-action on $V$. 
}
\end{rem}

\section{Orbifold $E$-functions}

\begin{dfn} 
{\rm Let $D= \sum_{j =1}^m d_j D_j$ be a $G$-invariant 
effective divisor on a smooth $G$-variety 
$V$ such that $(V,G)$ is a $G$-normal pair.  
Take  an arbitrary element  $g \in G$ and a connected component $W$ 
of $V^g$.  
Choose a point $x \in W$ and local $g$-invariant 
coordinates $z_1, \ldots, z_n$ at 
$x$ so that irreducible components of $Supp\, D$ containing 
$x$ are defined by local equations $z_i = 0$ for some 
$i \in \{1, \ldots, n\}$.  Let   $\delta_i $$( 1\leq i \leq n)$  
be the multiplicity 
of $D$ along $\{ z_i = 0 \}$ $( \{ \delta_1, \ldots, \delta_n \} 
\subset \{ 0, d_1, \ldots, d_m \} )$,  
and $e^{2\pi \sqrt{-1} \alpha_i}$ $( 1\leq i \leq n)$   the eigenvalue of 
the $g$-action  on $z_i$ ($\{ \alpha_1, \ldots, 
\alpha_n \} \subset  \Q \cap [0,1))$. 
We define the $D$-{\bf weight} of $g$ at  $W$ as 
\[ wt(g,W,D):= \sum_{i=1}^n \alpha_i(\delta_i +1). \]
If $D = 0$, then  
\[ wt(g,W):=  wt(g,W,0) = \sum_{i=1}^n \alpha_i \]
will be called simply the  {\bf weight} of $g$ at  $W$. 
Let $I^g$ be the subset of $g$-fixed elements in $I := \{1, \ldots, m \}$. 
For any subset $J \subset I^g$ we set 
\[ F(g,W,D^{\circ}_J; u,v) :=   
\prod_{j \in J} \frac{uv-1}{(uv)^{d_j+1} -1} 
E(W_J ; u,v), \]
where $W_J$ is the geometric quotient of $W \cap D_J^{\circ}$ modulo
the subgroup $C(g,W,J) \subset C(g)$ consisting of those elements in the 
centralizer of $g$ which leave the component $W \subset V^g$ and 
the subset $J \subset I^g$ invariant. 
} 
\end{dfn}

\begin{rem} 
{\rm
We note that  $wt(g,W,D)$ does not depend on the 
choice of a point $x \in W$. Moreover, if $h \in C(g)$ is 
an element in the centralizer of $g$ and $W' = hW$ is 
another  connected component of $V^g$, then  $wt(g,W',D) =  wt(g,W,D)$.  
} 
\end{rem}

\begin{dfn} 
{\rm 
We define the {\bf  orbifold $E$-function  of a  $G$-normal pair $(V, D)$} 
by the 
formula: 
\[ E_{\rm orb}(V,D,G; u,v) = \sum_{\{ g \}} 
\sum_{ \{ W \} } (uv)^{wt(g,W,D)} \sum_{J \subset I^g}  
F(g,W,D^{\circ}_J; u,v)
,  \]
where   $\{ g \}$ runs over all conjugacy classes in $G$, 
and $\{ W \}$ runs over the set of representatives of all 
$C(g)$-orbits in the set of connected components of $V^g$. 

In the case $D =0$, we call 
\[ E_{\rm orb}(V,G; u,v) : =  E_{\rm orb}(V,0, G; u,v) = 
\sum_{ \{ g \}   } 
\sum_{ \{W\} } (uv)^{wt(g,W_i)} E(W/C(g,W); u,v) , \]
the {\bf orbifold $E$-function} of a $G$-manifold $V$ (here 
$C(g,W)$ is the subgroup of all elements in $C(g)$ which leave 
the component $W \subset V^g$ invariant).  
}  
\end{dfn}

\begin{rem} 
{\rm Using the equalities  
\[ \frac{1}{|G|} \sum_{g \in G} \sum_{h \in C(g)}  
e(V^g \cap V^h) =  \sum_{\{ g \} \subset  G}   
\frac{1}{|C(g)|}   \sum_{h \in C(g)}     e(V^g \cap V^h)  = 
 \sum_{\{g \} \subset  G} e(V^g/C(g)), \]
one immediately 
obtains that $E_{\rm orb}(V, G; 1,1)$ equals 
the physicists' orbifold  Euler number $e(V,G)$ (see \ref{orb-ph}).  
\label{o-p}
} 
\end{rem} 

\begin{exam} 
{\rm Let $G := \mu_d$ a cyclic group of order $d$ acting by roots of 
unity on 
$V:= \C$. Then the corresponding orbifold  
$E$-function  
equals 
\[ E_{\rm orb}(V, G; u,v) = uv + \sum_{k =1}^{d-1} (uv)^{k/d} = 
(uv)^{1/d} + (uv)^{2/d} + \cdots + (uv)^{d-1/d} + uv. \]
} 
\label{mu-d1}
\end{exam}

\begin{lem} 
Let $V:= \C^r$ and $g \in {\rm GL}(r, \C)$ a linear authomorphism of finite
order. Denote by $V'$ the blow up of $V$ at $0$. Let $D \cong 
\P^{r-1}$ be the exceptional divisor in $V'$ and $\{ W_1, \ldots, W_s \}$ 
the set of connected components of $D^g$.  Then 
\[  \sum_{i =1}^s (uv)^{wt(g,W_i,D)} 
\frac{uv-1}{(uv)^r -1}
E(W_i; u,v) = (uv)^{wt(g,V^g)}. \]
\label{f1}
\end{lem}  

\noindent
{\em Proof.} Let  $\{ e^{2\pi \sqrt{-1} \alpha_i} \}$ 
$( 1\leq i \leq n)$ be the set of   the eigenvalues of $g$-action. 
Without loss of generality, we assume  
$0 \leq \alpha_1 \leq \cdots \leq \alpha_n <1$. We write the number $r$  
as a sum of $s$ positive integers $k_1 + \cdots + k_s$ where the numbers 
$k_1, \ldots, k_s$ are defined by the conditions 
$$\alpha_i = \alpha_{i+1} \; \Leftrightarrow\;  \exists j \in 
\{1, \ldots, s \}\; : \;  
k_1 + \cdots + k_j 
\leq i < k_1 + \cdots + k_j + k_{j+1} \]
and 
$$ \alpha_i < \alpha_{i+1}  \; \Leftrightarrow\;  \exists j \in 
\{1, \ldots, s \}\; : \;   
i+1 = k_1 + \cdots + k_j.$$ 
Then $D^g$ is a union of $s$ projectives linear subspaces $W_1, \ldots, 
W_s$,  where $W_j \cong  \P^{k_i -1}$ $(j \in \{1, \ldots, s \})$. 
By definition, we have $wt(g,V^g) = \sum_{i=1}^n \alpha_i$.  
By direct computations, one obtains   
$ wt(g, W_j,D) = k_1 + \cdots + k_{j-1} + \sum_{i=1}^r \alpha_i$. 
Hence, 
\[ \sum_{W_i \subset  {D}^g} (uv)^{wt(g,W_i,D)} E(W_i; u,v) = 
(uv)^{wt(g,V^g)} \sum_{j =1}^s (uv)^{ k_1 + \cdots + 
k_{j-1}}E(\P^{k_j -1}; u,v)  = \]
\[  (uv)^{wt(g,V^g)} \sum_{j =1}^s (uv)^{ k_1 + \cdots + k_{j-1}}(1 + 
(uv) + \cdots + (uv)^{k_j-1}) =  \]
\[ (uv)^{wt(g,V^g)} \sum_{l =0}^{r-1} (uv)^l = (uv)^{wt(g,V^g)} 
\frac{(uv)^r -1}{uv -1}. \] 
This completes the proof. \hfill $\Box$ 

\begin{lem} 
Let $V$ and $W$ be two smooth algebraic varieties having a regular 
action of a finite group $G$. Assume that $V$ is a Zariski locally trivial 
$\P^r$-bundle over $W$ such that the canonical 
projection $\pi\,: \,  V \to W$  is $G$-equivariant. Then 
\[ E(V/G; u,v) =  \frac{(uv)^r -1}{uv -1} E(W/G; u,v). \]
\label{f2}
\end{lem} 

\noindent
{\em Proof.} Let $H \subset G$ be a subgroup and $W(H):= 
\{ x \in W \, : \, St_G(x) = H \}$. Then $W \subset W$ is a locally closed 
subvariety, and $W$ admits a $G$-invariant stratification 
by locally closed strata 
\[ W = \bigcup_{ \{ H\} } W(\{H\}), \]
where $\{ H \}$ runs over the conjugacy classes of all subgroups in $G$ and 
$W(\{ H \}) := \bigcup_{ H' \in \{ H \} } W(H')$. Denote 
$V(\{ H \}) := \pi^{-1}(W(\{ H \}))$. Then  $V(\{ H \})$ is a $G$-equivariant 
$\P^r$-bundle over $W(\{H\})$ and we have isomorphisms 
 $V(\{ H \})/G \cong V(H)/N(H)$, 
 $W(\{ H \})/G \cong W(H)/N(H)$, where $W(H) := \pi^{-1}(W(H))$ and 
$N(H)$ is the normalizer of $H$ in $G$. Since $V(H)$ is a $N(H)$-equivariant 
$\P^r$-bundle over $W(H)$, it suffices to prove our statement for the 
case $G= N(H)$, $W = W(H)$, and  $V = V(H)$. Furthermore, we can restrict 
ourselves to the case when  $W$ is irreducible and $N(H)$ leaves $W$ 
invariant. The last conditions imply $N(H) = H$. Therefore, $W/G = W$ and 
the $H$-action on leaves each fiber of $\pi$ invariant. Hence, 
$E(V/G; u,v) = E(\P^r/H; u,v) E(W; u,v)$. Since all cohomology groups 
 of $\P^r$ have rank $1$ and they are generated by an effective algebraic 
cycle, we get $E(\P^r/H; u,v) =  E(\P^r; u,v)$. Thus, we have obtained 
the required formula for $E(\P^r/H; u,v)$.  \hfill $\Box$

\begin{theo} 
Let $(V,G)$ be a $G$-normal pair,  $Z \subset V$ a smooth $G$-invariant 
subvariety such that after the $G$-equivariant blow up $\psi \, : \, V' 
\rightarrow V$ with center in $Z$ one obtains a $G$-normal pair 
$(V', D')$, where $D'$ the effective divisor 
defined by the equality  
\[ K_{V'} = \psi^*(K_V - D) + D'. \]
Then
\[ E_{\rm orb}(V, D,G; u,v) =  E_{\rm orb}(V',D',G; u,v).  \]
\label{change}
\end{theo}

\noindent
{\em Proof.} Let $Z_1, \ldots, Z_k$ be the set of connected 
components of $Z$ and $D_1, \ldots, D_m$ the set of irreducible components of 
$Supp\, D$. Then $Supp\, D' = \psi^{-1} (Supp\, D) 
\cup D_{m+1} \cup \cdots \cup D_{m+k}$, where $D_{m+1}, \ldots, D_{m +k}$ are 
irreducible $\psi$-exceptional divisors over  
$Z_1, \ldots, Z_k$. It suffices to prove the equality 
\[ \;\; \forall g \in G: \;\;  \sum_{ \{ W \} } 
(uv)^{wt(g,W,D)} \sum_{J \subset I^g}  
F(g,W,D^{\circ}_J; u,v) =  \]
\[ = \sum_{ \{ W' \} } (uv)^{wt(g,W',D')} \sum_{J' \subset (I')^g}  
F(g,W',(D')^{\circ}_{J'}; u,v) ,    \]
where $I' = I \cup \{ m+1, \ldots, m+k \}$.  
We note that the $G$-equivariant 
mapping $\psi\, (V')^g \to V^g$ is surjective. Therefore, it 
suffices to prove the equality 
\[ (uv)^{wt(g,W,D)} \sum_{J \subset I^g}  
F(g,W,D^{\circ}_J; u,v) 
= \sum_{i =1}^l (uv)^{wt(g,W_i',D')} \sum_{J' \subset (I')^g}  
F(g,W_i',(D')^{\circ}_{J'}; u,v) ,    \]
where $W$ is a given  connected component of $V^g$ and 
$W_1', \ldots, W_l'$ are all connected components of $(V')^g$ such that 
$\psi(W_i') \subset W$ $( 1 \leq i \leq k)$. Since $\psi$ is an 
isomorphism over $W \setminus W \cap Z$ and the 
 $\psi$-exceptional divisors $D_{m+1}, \ldots, D_{m +k}$ 
are pairwise nonintersecting, it suffices to prove the equality 
\[ (uv)^{wt(g,W,D)}   
F(g,W \cap Z_j,D^{\circ}_J; u,v) 
= \sum_{i =1}^l (uv)^{wt(g,W_i',D')}   
F(g,W_i'\cap D_{m+j},(D')^{\circ}_{J'}; u,v),    \] 
where $j \in I$ and $J' = J \cup \{ j +m \}$. The last equality 
follows from Lemmas \ref{f1} and \ref{f2} using the fact that 
each  $W_i'\cap D_{m+j}$ is a locally trivial $\P^{k_i}$-bundle 
over $W \cap Z_j$. \hfill $\Box$

\section{Main theorems}

Let $V$ be a smooth $n$-dimensional 
algebraic variety, $G$ a finite group acting by 
regular authomorphism on $V$,  $X:=V/G$ it geometric quotient, and 
$\phi\, : \, V \rightarrow X$ the corresponding finite 
morphism. Then  $G$ acts on the set of irreducible components of 
the ramification divisor $\Lambda$ on $V$.  
Denote by  $\{ \Lambda_1, \ldots, \Lambda_k\}$ the set of representatives of
$G$-orbits  in the set of irreducible components of $Supp\, \Lambda$. 
Let $\nu_1 -1, \ldots, \nu_k-1$ 
be the multiplicities of  $\Lambda_1, \ldots, \Lambda_k$ in $\Lambda$ 
(the number $\nu_i$ equals the order of the cyclic 
intertia subgroup $St_G(\Lambda_i) \subset G$). 
Since $\phi\, : \, V \rightarrow X$ is a Galois covering, the multiplicity 
$\nu_i -1$ of $\Lambda_i$ depends only  on the $G$-orbit  
of $\Lambda_i$ in $Supp\, \Lambda$. 
We set  $\Delta_i := \phi(\Lambda_i)$ $( 1 \leq i \leq k)$   
and  consider the pair $(X, \Delta_X)$,   
where 
\[ \Delta_X  := \sum_{i =1}^k \left(\frac{\nu_i-1}{\nu_i}\right) \Delta_i 
\in Z_{n-1}(X) \otimes \Q. \]
By the ramification formula, we have   
\[ \phi^*(K_{X} + \Delta_X) = \phi^*K_{X} + \Lambda   = K_V.\]

\begin{prop} 
The pair  $(X, \Delta_X)$ is Kawamata log terminal. 
\end{prop} 

\noindent
{\em Proof.} 
 Let $\rho\, : \, Y \rightarrow X$ 
be a log resolution of singularities of $(X, \Delta_X)$
and 
\[ K_{Y} 
= \rho^* (K_{X} + \Delta_X) +  \sum_{i} a({D}_i,\Delta_X) {D}_i. \]
We consider the fiber product  $V_1 : = 
V \times_{X} Y$.  Then $V_1$ has  a natural finite Galois morphism  
$\phi_1\, :\, V_1 
\to Y$ and  a natural 
birational $G$-morphism  
$\rho_1\, : \, V_1 \to V$. We write  
\[ K_{V_1} = \rho_1^* K_V + \sum_{j=1}^m a( E_j,0) E_j,  \]
where $E_j$ runs over irreducible exceptional divisors of $\rho_1$.

By definition, the  multiplicity of any irreducible component $\Delta_i$ of 
$\Delta$ is equal to $(\nu_i -1)/\nu_i < 1$. Therefore, $a(D_i,\Delta_X) 
= - (\nu_j -1)/\nu_j  > -1$ if $\rho(D_i)$ 
coincides with an irreducible component $\Delta_j$ of  $Supp\, \Delta$.  
Now consider the case when  $\rho(D_i)$ is not an irreducible  
component of  $Supp\, \Delta$.
Denote by $E_j$ an irreducible divisor on  
$V_1$  such that $\phi_1(E_j) = D_i \subset Y$. 
Let $r_j$ be the ramification 
index of $\phi_1$ along $E_j$. By the ramification formula, one has 
$a(E_j,0)+1 = r_j( a(D_i, \Delta_X) +1)$. 
Since $V$ is smooth, we have $  a(E_j,0) \geq 1$ for all $j \in 
\{1, \ldots, m\}$. Therefore, $a(D_i, \Delta_X) = a(E_j,0)+1/r_j -1 > -1$.
\hfill $\Box$

\begin{dfn} 
{\rm Let  $V$ be a smooth algebraic variety having a regular action of 
a finite group $G$, and  $(X, \Delta_X)$ the pair  constructed above. 
Then  we call  $(X, \Delta_X)$ the 
{\bf Kawamata log terminal pair associated with } $(V, G)$. } 
\end{dfn}

\begin{exam} 
{\rm Let $G := \mu_d$ a cyclic group of order $d$ acting by roots of unity on 
$V:= \C$. Then $X=V/G \cong \C$ and $\Delta_X = \frac{d-1}{d} x_0$, where 
$x_0 \in X$ is the zero point. The stringy   
$E$-function of $(X, \Delta_X)$ 
equals 
\[ E_{\rm st}(X, \Delta_X; u,v) = (uv-1) + \frac{uv -1}{(uv)^{1/d} -1} = 
(uv)^{1/d} + (uv)^{2/d} + \cdots + (uv)^{d-1/d} + uv. \]
Thus, it coincides with the orbifold $E$-function from Example \ref{mu-d1}. 
 } 
\end{exam} 

Our next statements  show  the last  phenomenon in more general situations:

\begin{lem} 
Let $G \subset {\rm GL}(n, \C)$ be a finite abelian subgroup acting 
by diagonal matrices on $V:= \C^n$, and  $(X, \Delta)$  the Kawamata 
log terminal pair associated with $(V,G)$. 
Then 
\[  E_{\rm st}(X, \Delta_X; u,v) =  E_{\rm orb}(V, G; u,v). \]
\label{l-ab}
\end{lem} 

\noindent
{\em Proof.} 
First, we remark that the ramification locus $Supp\, \Lambda$ is contained 
in the union of the coordinate hyperplanes $\Lambda_i := 
\{ z_i = 0 \} \subset \C^n$ $(1 \leq i \leq n)$. 
Therefore, we can write $\Lambda = \sum_{i =1}^n \nu_i \Lambda_i$, where 
$\nu_i \geq 1$ $(1 \leq i \leq n)$. Second, we note that $X$ is a 
normal  affine toric variety corresponding 
to the cone $\sigma := \R_{\geq 0}^n$ and the lattice 
$$N := \Z^n + \sum_{g \in G} \Z(\alpha_1(g), \ldots, \alpha_n(g)),  $$
where $e^{2\pi \sqrt{-1} \alpha_1(g)}, \ldots, e^{2\pi \sqrt{-1} \alpha_n(g)}$ 
are the eigenvalues of $g$, $\{\alpha_1(g), \ldots, \alpha_n(g)\} \in 
\Q \cap [0,1)$. 
Moreover, $\Delta_X$ is a torus invariant divisor on $X$.
Let us denote by $\{ e_1, \ldots, e_n \}$ the standard basis of $\Z^n$. 
Then   the $\Q$-divisor $K_X + \Delta_X$ corresponds to a linear 
function $\varphi_{K, \Delta}$ which has value $1$ on each 
$e_i$ $(1 \leq i \leq n)$. By \ref{tor-klt1}, $(X, \Delta_X)$ is a torus 
invariant Kawamata log terminal pair. By \ref{st-tor}, we obtain 
\[   E_{\rm st}(X, \Delta_X; u,v)= (uv-1)^n \sum_{p \in N \cap \sigma} 
(uv)^{-\varphi_{K,\Delta}}. \] 
We set  $f_i := (1/\nu_i)e_i$ $(1 \leq i \leq n)$. 
Then the system of vectors 
$\{ f_1, \ldots, f_n \} \subset  N$ 
generates a sublattice $N' \subset N$ containing  
$\Z^n$. Denote by ${\cal R}:= \{ v_1, \ldots, v_r \} \subset  N$ the  set   
of representatives of $N/N'$, where each element $v \in 
{\cal R}$ has  a form $v= 
\sum_{ i=1}^n \lambda_i(v)  f_i $ $(0 \leq \lambda_i < 1)$. 
Then, by summing a multidimensional geometric series, we obtain  
\[  (uv-1)^n \sum_{p \in (v + N' )\cap \sigma} 
(uv)^{-\varphi_{K,\Delta}}  =  (uv-1)^n \left( 
(uv)^{-\sum_{ i=1}^n \lambda_i(v)/\nu_i} \right) \prod_{i=1}^n 
\frac{1}{1 - (uv)^{-1/\nu_i}} \]
(we used the property $\varphi_{K, \Delta}(f_i) = 1/\nu_i$,  
$1 \leq i \leq n$).     
Thus,  we have 
\[  E_{\rm st}(X, \Delta_X; u,v) = \left( \sum_{ v \in {\cal R}} 
(uv)^{- \sum_{ i=1}^n \lambda_i(v))/\nu_i} \right) \prod_{i=1}^n 
\frac{(uv-1)}{1 - (uv)^{-1/\nu_i}} =  \]
\[ (uv)^n  \left( \sum_{ v \in {\cal R}} 
(uv)^{- \sum_{ i=1}^n \lambda_i(v))/\nu_i} \right) 
\prod_{i=1}^n \left(1 + (uv)^{-1/\nu_i} + \cdots + (uv)^{-(\nu-1)/\nu_i}   
\right) =  \]    
\[ =  (uv)^n \sum_{ g \in G} (uv)^{ - \sum_{i =1}^n \alpha_i(g)} = 
E_{\rm orb}(V, G; u,v). \]
\hfill $\Box$

Now we come to our main theorem: 

\begin{theo}
Let $G$ be a finite group acting regularly on a smooth algebraic 
variety $V$ and  
$(X, \Delta)$  the Kawamata 
log terminal pair associated with $(V,G)$. 
Then  
\[   E_{\rm st}(X, \Delta_X; u,v) =  E_{\rm orb}(V, G; u,v).  
\]
\label{equal3}
\end{theo}
 
\noindent 
{\em Proof.} Let $(V^{ab}, D)$ be the canonical abelianization 
of the $G$-normal pair $(V,0)$, $D = K_{V^{ab}} - \psi^* K_V = 
\sum_{i =1}^m d_i D_i$. Denote by $\phi^{ab}$ the finite morphism 
$V^{ab} \to Y: = V^{ab}/G$. Then $\psi$ induces a  
birational proper morphism  $\overline{\psi}\, : \, Y \to X$ which 
can be considered as a partial desingularization of $X$. Let $W_1, 
\ldots, W_l$ be representatives of $G$-orbits in the set 
$\{ D_1, \ldots, D_m \}$ 
$( l \leq m)$, and $\overline{W}_1, \ldots, \overline{W}_l$ 
their $\phi^{ab}$-images in $Y$. By the ramification formula, we have 
\[ K_{Y} 
= (\overline{\psi})^* (K_{X} + \Delta_X) + 
\sum_{j=1}^{l}  \left( \frac{d_j+1}{r_j} -1 \right) \overline{W}_j
+  \sum_{i =l+1}^{l+k}  
\left(\frac{1}{\nu_i} -1\right)\overline{W_i} ,   \]
where $\overline{W_{l+i}}$ is the $\phi^{ab}$-image of 
$\psi^{-1}(\Lambda_i) \subset V^{ab}$ in $Y$, and 
 $r_j$ is the  order of the ramification of $W_j$ over  
$\overline{W}_j$. We set $I_1 := \{1, \ldots,  \}$,      
$I_2 := \{l+1, \ldots, l+k \}$ and $I := I_1 \cup I_2$.  
For any subset $J \subset I$ we set $J_1 = I_1 \cap J$ and  
$J_2 = I_2 \cap J$. Denote by   $G(J)$ the $G$-stabilizer
of a point $x \in V^{ab}$ such that $\phi^{ab}(x) 
\in \overline{W}_J^{\circ}$ and set 
\[ S(J; u,v) :=  \sum_{ g \in G(J) } (uv)^{wt(g,x,D)}. \] 
It is easy to see that if $x' \in V^{ab}$ is another 
point such that  $\phi^{ab}(x') \in \overline{W}_J^{\circ}$, then 
$St_G(x')$ is conjugate to  $St_G(x)$, i.e., $G(J)$ depends only on 
$J$, but not on the choice
of a point $x \in (\phi^{ab})^{-1}(\overline{W}_J^{\circ})$. 
Let    $G'(J)$ be the subgroup in $G(J)$ 
generated by the cyclic inertia subgroups $St_G({W}_j)$ $(j \in J_1)$
and  $St_G(\psi^{-1}({\Lambda}_{j-l}))$ $(j \in J_2)$; i.e., we have 
$G'(J) \cong \prod_{j \in J_1} \mu_{r_j}  \prod_{j \in J_2} \mu_{\nu_j}, 
$ and 
\begin{equation} 
S'(J; u,v):= \sum_{ g \in G'(J) } (uv)^{wt(g,x,D)} =  
\prod_{j \in J_1} \frac{ (uv)^{d_j+1} -1}{ (uv)^{(d_j+1)/r_j} -1}  
\prod_{j \in J_2}  \frac{ uv -1}{ (uv)^{1/\nu_j}-1 } 
\label{e1}
\end{equation} 
By \ref{change}, we have 
\[ E_{\rm orb}(V, G; u,v) =  E_{\rm orb}(V^{ab},D, G; u,v) = 
\sum_{J \subset I} S(J; u,v)
\prod_{j \in J_1} 
\frac{ uv -1}{ (uv)^{d_j+1} -1}E(\overline{W}_J^{\circ}; u,v).   \]
Since the singularities along $\overline{W}_J^{\circ}$ 
are toroidal (cf. \cite{D}), in follows from 
\ref{l-ab} that  
\[  E_{\rm st}(X, \Delta_X; u,v) = \sum_{J \subset I} 
\overline{S}(J; u,v)  
\prod_{j \in J_1}  \frac{ uv -1}{ (uv)^{(d_j+1)/r_j} -1} 
\prod_{j \in J_2}  \frac{ uv -1}{ (uv)^{1/\nu_j}-1 }
E(\overline{W}_J^{\circ}; u,v),   
\] 
where $\overline{S}(J; u,v)S'(J; u,v)= S(J; u,v)$. It remains to apply 
(\ref{e1}).   
\hfill $\Box$

\bigskip

\noindent
{\bf  Proof of Theorem \ref{equals}:} 
The statement immediately follows from \ref{o-p} 
and \ref{equal3} by taking limits: 
\[ e_{\rm st}(X, \Delta_X) = 
\lim_{u,v \to 1} E_{\rm st}(X, \Delta_X; u,v) = 
\lim_{u,v \to 1} E_{\rm orb}(V, G; u,v) = e(V,G). \]
\hfill $\Box$

\begin{coro} 
Let $X$ be a normal complex algebraic surface with at worst log terminal 
singularities. Then 
\[ e_{\rm st}(X) = e(X \setminus X_{sing} ) + \sum_{x \in  
X_{sing}} c_x, \]
where $c_x$ is the number of conjugacy classes in the 
local fundamental group of $X \setminus \{ x \}$. In particular, 
$e_{\rm st}(X)$ is always an integer.  
\end{coro} 

\noindent
{\em Proof.} It is well-known that a germ of a singular point 
$x \in X_{sing}$ is isomorphic to a germ of $0$ in 
$\C^2/G_x$ where $G_x \subset {\rm GL}(2, \C)$ is a finite subgroup
( $G_x$ is isomorphic to the local fundamental group of $X \setminus 
\{x\}$).
Therefore, we have  $J_{\infty}(X,x) \cong 
J_{\infty}(\C^2/G_p, 0)$. Let $\rho\, : \, Y \to X$ be a  
resolution of singularities, $D_1, \ldots, D_m$ are exceptional divisors
over $x \in X$, $\{ a_1, \ldots, a_m \}$ their discrepancies, and 
$I = \{ 1, \ldots, m \}$. By \ref{indep} and \ref{equals}, the number 
\[ e_{\rm st}(x) := \sum_{J \subset I} e(D^{\circ}_J) \prod_{j \in J} 
\frac{1}{a_j +1} \]
does not depend on the choice of a resolution and equals $c_x$. 
\hfill $\Box$

\section{Cohomological McKay correspondence}

\begin{dfn} 
{\rm Let $G \subset {\rm SL}(n, \C)$ be a finite subgroup acting linearly 
on $V:= \C^n$ and $X:= V/G$. A resolution of singularities 
$\rho\, : \, Y \to X$ is called {\bf crepant}  if the canonical class $K_Y$ 
is trivial.}   
\end{dfn}

\begin{prop}
Let  ${\C}^* \times X \rightarrow X$ be the regular ${\C}^*$-action on 
$X$ induced by the action of scalar matrices  on ${\C}^n$. Assume that 
there exists a crepant resolution of singularities  
$\rho\, : \, Y \to X$.  Then the  
${\C}^*$-action on $X$ extends uniquely 
to a regular ${\C}^*$ -action on  $Y$. 
\end{prop}

\noindent
{\em Proof.}   
Since $Y$ is birational to $X$, the ${\C}^*$-action on $X$ extends 
uniquely to a rational ${\C}^*$-action ${\C}^* \times Y \dasharrow  Y$. 
It remains to show that it  is  regular. Let $\{ D_1, \ldots,
 D_m \}$  be the set of all irreducible 
divisors on $Y$ in the exceptional locus of $\rho$. 
It was shown in \cite{IR} that the corresponding 
discrete valuations ${\cal V}_{D_1}, 
\ldots,{\cal V}_{D_m}$ of the field of rational functions on $Y$ are 
determined uniquely. Since the algebraic group $\C^*$ is connected, every 
such a  
valuation   ${\cal V}_{D_1}, \ldots,{\cal V}_{D_m}$ must be  
invariant under the 
rational ${\C}^*$-action on $Y$. Therefore, the rational 
${\C}^*$-action on $Y$ can be extended to a regular action on some 
Zariski dense open subsets $U_j \subset D_j$ $( j=1, \ldots, m)$, i.e., 
the rational  ${\C}^*$-action on $Y$  is regular outside some 
Zariski closed subset 
$$Z := \bigcup_{j =1}^m (D_j \setminus U_j) \subset Y,\; \;\;  
codim_Y\, Z \geq 2.$$  
Let $TY$ be 
the tangent vector bundle over $Y$. By the extension theorem of   
Hartogs, the restriction mapping  on 
global sections  
$\Gamma(Y, TY) \rightarrow \Gamma(Y \setminus Z, TY)$ is 
bijective. Hence,  the regular 
vector field $\eta  \in \Gamma(Y \setminus Z, TY)$ 
corresponding to the regular 
${\C}^*$-action on $Y  \setminus  Z$  
extends to a regular vector field on the whole variety $Y$. The latter 
shows that the  
${\C}^*$-action on $Y \setminus Z$ extends to a regular action 
on the whole $Y$. \hfill $\Box$

\begin{lem}
Let $V$ be a smooth algebraic variety,  and $W = \bigcup_j W_j$
a stratification of $W$ by locally closed irreducible subvarieties. 
Assume that the Hodge structure in the cohomology with 
compact supports  $H^{i}_c(W_j, \Q)$ is pure for all $i,j$. 
Then the Hodge structure in  $H^{i}_c(W, \Q)$
is  pure for all $i$. 
\label{decomp}
\end{lem}

\noindent 
{\em Proof.}  The statement  
follows by induction using tha fact that for any closed subvariety 
$W' \subset W$ the 
long exact cohomology sequence   
\[ \rightarrow H^{i-1}_c (W')   \rightarrow H^{i}_c (W \setminus 
W')   
\rightarrow H^{i}_c (W)  
  \rightarrow H^{i}_c (W')  \rightarrow H^{i+1}_c (W \setminus 
W') \rightarrow \]
respects the Hodge structure. 
\hfill $\Box$  
\bigskip

The following statement was conjectured in \cite{BD} (see also 
\cite{IR}): 

\begin{theo}
Let $G \subset {\rm SL}(n, {\C})$ be a finite 
subgroup. Assume that there exists  a crepant desingularization 
$\rho \;:\;Y \rightarrow X := {\C}^n/G$. 
Then the Hodge structure in the cohomology  
$H^*(Y, {\C})$ is pure. Moreover, 
$H^{2i +1}(Y, {\C}) = 0$,  
 $H^{2i}(Y, {\C})$ has the Hodge type $(i,i)$ for all $i$, and 
the dimension of $H^{2i}(Y, {\C})$ is equal to 
the number of conjugacy classes $\{ g \} \subset G$ having the  
weight $wt(g)=i$. 
\label{mck}
\end{theo}

\noindent 
{\em Proof.} 
Let $Y^{{\C}^*}$ be the fixed point set 
of the ${\C}^*$-action on $Y$, 
$Y^{{\C}^*} = \bigcup_{ j =1}^l Y_j$
a decomposition of $Y^{{\C}^*}$ in its connected components, 
$Y_0 := \rho^{-1}(x_0) \subset X$, where $x_0$ is the image of $0 \in \C^n$ 
modulo $G$. Since $Y_0$ is 
the fiber over the unique ${\C}^*$-fixed 
point $x_0 \in X$, we have  $Y^{{\C}^*} \subset Y_0$. Therefore 
 $Y^{{\C}^*}$ is compact. Since the fixed point subvariety 
$Y^{{\C}^*}$ is smooth and compact, the cohomology of  
every  connected  component $Y_1, \ldots, Y_k$ of   $Y^{{\C}^*}$ 
have pure Hodge structure. 
Consider the Bialynicki-Birula cellular decomposition \cite{B-B}: 
$Y = \bigcup_{ j =1}^l W_j$, where  
$W_j = \{ y \in Y \; :\; \lim_{z \rightarrow 0} z(y) \in Y_j, \; z 
\in {\C}^* \}$.
Since every $W_j$ is a vector bundle over $Y_j$, the   groups  
$H^i_c(W_j, \C)$ have   pure Hodge structures for all  $i,j$.  
By \ref{decomp}, the Hodge structure 
in $H^i_c(Y, \C)$ is pure for all $i$. 

Denote by  $C_i(G)$ the number of 
conjugacy classes $\{ g \} \subset G$ having the  
weight $wt(g)=i$. Since $G$ is contained in 
${\rm SL}(n, \C)$, the ramification divisors $\Lambda$
and  $\Delta_X$ are zero. By \ref{int-for} and 
\ref{equal3}, we have 
\[ E(Y; u,v) = E_{\rm st}(X, 0; u,v) =  E_{\rm orb}(\C^n, G; u,v). \]
Using the purity of  $H^i_c(Y, \C)$ and the fact that the 
Poincar{\'e} duality 
\[  H^{2n-i}_c(Y, \C) \otimes  H^{i}(Y, \C) \to  H^{2n}_c(Y, \C) 
\cong \C(n) \]
respects the Hodge structure, it remains to show that 
\begin{equation}
E_{\rm orb}(\C^n, G; u,v) =  \sum_{\{ g \}} C_i(G)  (uv)^{n-i}.
\label{for}
\end{equation} 
Indeed, we have 
$E_{\rm orb}(V, G; u,v) = \sum_{\{g \}} (uv)^{wt(g,V^g)} E(V^g/C(g); u,v)$, 
where $V:= \C^n$. Since $V^g$ is a linear subspace of dimension 
$k(g):= dim\, Ker(g - id)$, we obtain $E(V^g/C(g); u,v) = (uv)^{k(g)}$. 
Hence, 
$$(uv)^{wt(g,V^g)} E(V^g/C(g); u,v) = (uv)^{n - wt(g^{-1}, V^g)}. $$
The summing over all  conjugacy classes $\{ g^{-1} \}$ implies 
(\ref{for}).  
\hfill $\Box$ 
\bigskip
\bigskip

\noindent
{\bf  Proof of Theorem \ref{equals2}:} Now it  follows immediatelly from 
\ref{mck}.  \hfill $\Box$

\end{document}